\documentclass[12pt]{amsart}
\usepackage{amsfonts,amsmath,amssymb,amsthm,verbatim}
\usepackage[latin1]{inputenc}
\theoremstyle{plain} %% This is the default
\newtheorem{thm}{Theorem}%[section]

\newtheorem{lem}[thm]{Lemma}

\theoremstyle{definition}

\newtheorem{defn}{Definition}

\theoremstyle{remark}

\newcommand{\thmref}[1]{Theorem~\ref{#1}}

\newcommand{\lemref}[1]{Lemma~\ref{#1}}

\newcommand{\x}{\mathbf{x}}

\ifx\undefined\xspace \usepackage{xspace} \fi
\newcommand{\qr}[1]{\eqref{#1}}

 \newcommand{\RR}{{\mathbb R}}
\newcommand{\ZZ}{{\mathbb Z}}

\newcommand{\ifemptythenelse}[3]{%
\def\epty{} %
\def\atempempty{#1} %
\ifx\atempempty\epty %
#2 %
\else %
#3 %
\fi
}

\renewcommand{\d}{\,{d}}

\numberwithin{equation}{section}
\numberwithin{thm}{section}
\newsavebox{\Chai}
\sbox{\Chai}{\raisebox{0.3ex}[1.7ex][0.75ex]{{$\chi$}}}

\renewcommand{\L}[1]{L_{#1}}

\renewcommand{\x}{x}
\newcommand{\CF}{\mathcal{F}}

\newcommand{\var}{\operatorname{var}}

\let\tl=\tilde
\let\mc=\mathcal
\let\X=X
\let\x=x
\let\T=T

\newlength{\ppph}
\def\ppp#1\par{\par\paragraph{#1}}
\begin{document}
\title[Square summability of variations and the transfer operator]{Square summability of variations and convergence of the transfer operator}
\author{Anders Johansson and Anders \"Oberg}

\address{Anders Johansson\\ Division of Mathematics and Statistics\\
  University of G\"avle\\ SE-801 76 G\"avle\\ Sweden} \email{ajj@hig.se}

\address{Anders \"Oberg\\
  Department of Mathematics\\ Uppsala University\\ P.O. Box 480 \\ SE-751 06
  Uppsala \\ Sweden} \email{anders@math.uu.se}


\maketitle
\begin{abstract}
  In this paper we study the one-sided shift operator on a state space
  defined by a finite alphabet. Using a scheme developed by Walters
  \cite{walters2}, we prove that the sequence of iterates of the transfer operator converges under
  {\em square summability of variations} of the $g$-function, a
  condition which gave uniqueness of a $g$-measure in \cite{johob}. We
  also prove uniqueness of so-called $G$-measures, introduced by Brown
  and Dooley \cite{brown}, under square summability of variations.
\end{abstract}

\section{Introduction}\noindent
We consider the left one-sided shift map $T$ on the state space
$X^+=S^{\mathbb Z_+}$, where $S$ is a finite set. Thus $T$ acts on
elements $x$ of $X^+$, $x=(x_0,x_1, x_2,\ldots)$, in the following way
(each $x_i$ belongs to $S$):
$$T(x_0, x_1, x_2,\ldots)=(x_1, x_2,\ldots).$$

A $g$-measure associated to a continuous function $g: S^{\mathbb
  Z_+}\to [0,1]$ is a $\T$-invariant measure $\mu$ in the space of
Borel probability measures $\mc PX^+$, such that $g = {d\mu}/{d\mu
  \circ T}$. Such a function is referred to as a $g$-function and
satisfies $\sum_{y\in \T^{-1}x} g(y) = 1$, for all $x\in X^+$, and
thus can be viewed as a probability transition function for the local
inverses of $T$. Keane \cite{keane} introduced the concept of
$g$-measures into Ergodic Theory.

Existence of at least one $g$-measure is automatic if $S$ is finite
and if $g>0$, since then $S^{\mathbb Z_+}$ is compact, and the
continuity of $g$ implies that $\inf g>0$. If $S$ is countably
infinite, then existence is no longer automatic, but a weak sufficient
condition was given in \cite{jopp}.

Recent results concerning necessary conditions for uniqueness of a
$g$-measure have been given by Berger, Hoffman and Sidoravicius
\cite{berger} and also by Hulse \cite{hulse}.

Necessary and sufficient conditions are sometimes, but not always,
centered around conditions on the variations of the $g$-function on
small parts of its domain. For $f: X^+ \to \RR$, the variation is
defined as
$$
\var_{n} f = \sup_{x\sim_n y} |f(x) - f(y)|,
$$ 
where $x\sim_n y$ if $x,y\in X^+$ coincide in the first $n$
coordinates.

Walters \cite{walters1} proved uniqueness and rates of convergence
under summable variations of $g$-functions, or rather their
logarithms, $\log g$, which amounts to the same thing, if $\inf g >0$.

When $S$ is finite, a sufficient condition for uniqueness of a
$g$-measure, {\em square summability of variations} of the
$g$-function, was found by the present authors \cite{johob} and this
result was extended to {\em countable state shifts} in Johansson,
\"Oberg and Pollicott \cite{jopp}.

\newcommand{\svar}{\operatorname{svar}} \newcommand{\Ti}[1]{T^{-#1}}  
Square summability of variations of a $g$-function then means that
\begin{equation}
    \sum_{n=1}^{\infty} \left(\var_{n} f\right)^2<\infty.
\end{equation}

In this paper we prove convergence of the iterates of the
transfer (transition) operator, which is defined by
\begin{equation}
    \mathcal\L f(x)=\sum_{y\in T^{-1}x} g(y)f(y),
\end{equation}
for continuous functions $f$ on $X^+$.

More specifically, we prove that with finitely many symbols in $S$ and
under square summability of variations of a strictly positive
$g$-function, we have $\mathcal\L^n f(x)\to \int f\, d\mu$, uniformly
in $x$.

To accomplish this, we use a theory developed by Walters in
\cite{walters2} to study the sequence of iterates of an operator $P$,
which is defined as
\begin{equation}
    P^n f(x)=\left(\mathcal\L^n f\right)(T^n x),
\end{equation}
for $n\geq 1$.  

Walters introduces an adjoint sequence of operators
$P^{(n)^{\ast}}:\mc PX^+\to \mc PX^+$, $n\geq 1$. Define
$K^n_g=\{\nu\in \mc PX^+|P^{(n)^{\ast}}\nu=\nu\}$. Then (Theorem
2.1(v) of \cite{walters2}) implies that $K^n_g=P^{(n)^{\ast}}\mc PX^+$
and that $K^1_g\supset K^2_g\supset \ldots$. Walters defined
$K_g=\bigcap_{n=1}^{\infty}K^n_g$, which is the set of $G$-measures in
the terminology of Brown and Dooley \cite{brown}. See also Fan
\cite{fan} and Dooley and Stenflo \cite{dooley} for some recent
contributions to uniqueness of $G$-measures. The motivation for
studying $G$-measures has to some extent been to obtain a more general
framework for studying Riesz products, see \cite{fan}, for example.

Our strategy is now to prove that there exists a unique measure in
$K_g$ (it is trivial that any $g$-measure must be an element of
$K_g$). Then we use Theorem 2.9 of \cite{walters2} to conclude that
$\mathcal\L^n f(x)\to \int f\, d\mu$, where $\mu$ is the unique
element of $K_g$. This follows by a compactness argument, so in order
to prove the convergence of $\mathcal\L^n f$, we have assumed that the
symbol space $S$ is finite.

We need to define a probability measure, which is a Markov chain on $X^+$, which we call a $g$-chain, in order to reverse the dynamics in the absence of stationarity; it is not known if (when) the measures in $K_g$ are invariant under $T$.

To conclude the proof, we use the same method as in \cite{johob} to prove that two extremal measures in $K_g$ must be absolutely continuous, and hence that there exists only one such measure. This method was developed further in \cite{jopp} to cover the case of countable state shifts (thus improving on the results of \cite{johob}) for $g$-measures, using more sophisticated methods to prove absolute continuity of measures. These methods had appeared much earlier in the literature by Shiryaev and co-authors, see e.g.\ the survey by Shiryaev \cite{shiryaev}.

We end the paper with two open questions.

\section{Results and proofs}\noindent

\begin{thm}\label{G}
  Assume we have a finite symbol set $S$ and let $g>0$ satisfy
$$\sum_{n=1}^\infty \left(\var_n \log g\right)^2 
   \asymp \sum_{n=1}^\infty \left(\var_n g\right)^2<\infty$$
and let $\mu$ be the corresponding unique $g$-measure. Then
$K_{g}=\{\mu\}$, or equivalently: there exists a unique $G$-measure,
which is the same as the unique $g$-measure.
\end{thm}

The following theorem is in view of \cite{walters2} a corollary
of \thmref{G}. It uses the connection of Theorem 2.9 of \cite{walters2} and we copy part of its
proof for the convenience of the reader, since it is here explained
why we have assumed that we have finitely many symbols
(one has to have a compact state space, or at least some tightness property)
so that there exists a certain convergent subsequence for a sequence
of measures.

\begin{thm}\label{conv}
  Assume we have a finite symbol set $S$ and let $g>0$ be a
  $g$-function with square summable variation and let $\mu$ be the
  corresponding unique $g$-measure. Then $\sup_x |\mathcal\L^n
  f(x)-\int f \d\mu|\to 0$.
\end{thm}

\begin{proof} (Borrowed from Theorem 2.9 of \cite{walters2}.)  Under
  the same assumptions, \thmref{G} gave us a unique measure in $K_g$.
  From this it follows that $\sup_x|P^nf(x)-\int f\d\mu |\to 0$, since
  otherwise, there exists a continuous function $f$ and an
  $\epsilon$ and sequences $\{n_k\}$ and $\{x_k\}$ such that
$$
\left|\int f\d\left(P^{(n_k)^\ast}\delta_{x_k}\right)-\int
  f\d\mu\right| \geq \epsilon \; \; \forall k\geq 1.
$$ 
We can now, due to compactness, pick a subsequence $\{k_j \}$ of
$\{k\}$ such that the sequence of measures
$\{P^{(n_{k_j})^\ast}\delta_{x_{k_j}}\}$ converges to $\nu\in \mc
PX^+$. But then since $K_g^{(n)}=P^{(n)^\ast} \mc PX^+$, we also
have $\nu \in K_g$, a contradiction. Since
$P^{(n)}f(x)=\left(\mathcal\L^n f\right)\left(T^n x\right)$, we also
have $\sup_x|\mathcal\L^nf(x)-\int f\d\mu |\to 0$.
\end{proof}

\subsection{Proof of  \thmref{G}}
We begin with some preliminary terminology; see \cite{jopp} for a more
thorough exposition. For a pair of probability measures $\nu,\tl\nu
\in\mc P\X$ and some filtration, $\{\CF_n\}$, let $Z_n(x) = Z_n(x;
\tl\nu,\nu,\CF_n)$ be the likelihood ratio martingale
\[ Z_n(x;\tl\nu,\nu,\{\CF_n\}) 
= \dfrac{\d\tl\nu\vert_{\CF_{n}}}{\d\nu\vert_{\CF_{n}}}(x), \]
on $\CF_n$, where we assume that $\tl\nu$ and $\nu$ are locally absolute continuous on $\CF_n$ (see \cite{jopp}). 

Now, let $\X=S^{\ZZ}$ and $X^+=S^{{\mathbb Z}_+}$ and extend the
one-sided shift $T$ on $X^+$ to a two-sided shift on $\X$. For $a,b
\in \ZZ$ let $\Pi_{a,b}(x)$ be the mapping taking $x\in\X$ to
$(x_a,x_{a+1},\dots,x_b)\in S^{b-a}$ and define the natural projection
$\Pi_+: \X\to \X^+$ taking bi-infinite sequences $x$ in $\X$ to
one-sided ones in $\X^+$ by \(\displaystyle \Pi_+\left(
  (x_i)_{i=-\infty}^{i=\infty} \right) = (x_0,x_1,\dots). \) Let also
$\CF_a^b$ be the algebra generated by $\Pi_{a,b}$ and $\CF_n^-:=
\CF_{-n}^{-1}$ (the ``forward'' algebra) and $\CF_n^+:=\CF_{0}^{n-1}$
(the ``backward'' algebra). A {\em cylinder set} is a set of the form $[x]_a^b=\Pi_{a,b}^{-1}\Pi_{a,b}(x)$.

We will now define a certain Markov chain, a {\em $g$-chain} in order to be able to go forward in time, since we may not assume that the measures in $K_g$ defined by Walters are invariant under $T$. 

\begin{defn}
A $g$-chain on $\X=S^{\mathbb Z}$ is a probability measure $\nu \in \mc P\X$
such that, for all $n$ in $\ZZ$,
\begin{equation}\label{gchain}
  \nu(x_n | x_{n+1},x_{n+2},\dots) = g(x_n, x_{n+1},x_{n+2},\dots ).
\end{equation}
A forward $g$-chain is a probability measure $\nu$ on $\X$ satisfying \eqref{gchain} for $n\leq -1$.
\qed
\end{defn}

For a $g$-chain $\nu$, the distribution under $\nu$ of the process $x^{(t)} \in \X^+$,
$t\in\ZZ$, is defined by
\[ x^{(t)} := \Pi_+(T^{-t}x)=(x_{-t},x_{-t+1},\ldots ), \] 
is that of a Markov chain such that the transition probabilities are
given by $g$ (and the transition operator of the chain is $\mc L$). That is, for all $t\in\ZZ$, 
$$\nu(x^{(t)}|x^{(t-1)}) = g(x^{(t)}).$$ 
The same holds for $t\geq 1$ in the case of a forward $g$-chain
$\nu$.

\begin{lem}\label{corr}
A probability measure $\nu\in\mc P\X$ is a $g$-chain only if $\nu\circ
\Pi_+^{-1} \in \mc P\X^+$ is an element of the set $K_g$ of
eigenmeasures as defined by Walters in \cite{walters2}. And vice
versa, any $\nu$ in $K_g$ corresponds, by extension, to a
unique $g$-chain.
\end{lem}
\begin{proof}
Suppose that we have $P^{\ast}\nu=\nu$. This is equivalent with $\mc L^{\ast}(\nu \circ T^{-1})=\nu$. By interpreting the conditional probability
$\nu(x_0|x_1,x_2,\ldots)$ as the transition probability $\nu(x^{(0)}|x^{(-1)})$ from the Markov chain defined above, we reach the conclusion that $$\nu(x_0|x_1,x_2,\ldots)=g(x_0,x_1,\ldots).$$
Thus, if $P^{(n)^{\ast}}\nu=\nu$, we have for $0\leq k\leq n-1$
\begin{equation}\label{chain}
\nu(x_k | x_{k+1},x_{k+2},\dots) = g(x_k, x_{k+1},x_{k+2},\dots ).
\end{equation}
Since $K^n_g=P^{(n)^{\ast}}\mc PX^+$ and $K_g=\bigcap_{n=1}^{\infty}K^n_g$, we have precisely \eqref{chain} for all $k\geq 0$. Using measures in the nonempty set $K_g$ as initial distributions, we may define the full $g$-chain, by making a unique canonical extension (Neveu \cite{neveu}, p.\ 83, the corollary), so that \eqref{chain} is true for $k\leq -1$. Hence \eqref{chain} holds for all $k\in {\mathbb Z}$.
\end{proof} 

From now on, we do not distinguish between a
$g$-chain $\nu$ and its one-sided restriction $\nu\circ \Pi_+^{-1}$.

\begin{lem}\label{extremal}
  The set $K_g$ is a non-empty convex subset with
  mutually singular extreme points.
\end{lem}
\begin{proof}
See Theorem 2.11 and its proof in Walters \cite{walters2}.
\end{proof}

We work under the assumption that $g>0$ and hence that any two
$g$-chains be will locally a.c.\ on any of the algebras $\CF_a^b$, where $a$ and $b$ are finite. Given two
$g$-chains, $\nu,\tl\nu\in K_g$, let
\begin{equation}\label{xidef}
  \xi_n(x) := Z_n(x; \tl\nu,\nu,\{\CF_n^+\})
\end{equation}
Note that $\xi_n$ is the likelihood-ratio martingale and since it is a
positive $L_1(\nu)$-bounded $\nu$-martingale, we know that it
converges $\nu$-almost surely. We want to show that it is uniformly
integrable, {\em UI}, with respect to $\nu$, which gives an $L_1$- density between $\tl\nu$ and $\nu$ on $\CF^+=\lim \CF_n^+$.  And thus, according to \lemref{extremal}, a contradiction if $\nu$ and $\tl\nu$ are chosen to be distinct and extremal in $K_g$.

It was shown in \cite{johob} (p.\ 595) that the {\em UI} property of a
likelihood-ratio martingale such as $\xi_n$ amounts to show that
\begin{equation}
  \lim_{K\to\infty} \sup_n \tl\nu( \log \xi_n > K ) = 0.\label{tight}
\end{equation}

To see \qr{tight}, we note that, for a fixed $m$, a translation $m$ steps to the left of both the point and the measure gives that
\begin{equation}\label{xidef}
  \xi_m(x) = Z_m(T^m x; \tl\nu\circ T^{-m},\nu\circ T^{-m},\{\CF_m^-\}).
\end{equation}
and hence the law of $\xi_m$ under $\tl\nu$,
$\tl\nu\circ\xi_m^{-1}$, equals the distribution
$(\tl\nu\circ T^{-m}) \circ \zeta_{m,m}^{-1}$ 
where $(\zeta_{m,n}(x):n\in\ZZ_+)$ is the forward likelihood ratio
martingale ($m$ is fixed)
$$
  \zeta_{m,n}(x)= 
  Z_n(\, x ;\,\tl\nu\circ T^{-m},\nu\circ T^{-m},\{\CF_n^-\}).  
$$
This means that we start with two extremal measures $\nu, \tl\nu \in K_g$ translated $m$ times to the left, i.e., $\nu\circ T^{-m}$ and 
$\tl\nu \circ T^{-m}$, which can be extended to well-defined $g$-chains and then we may go forward along $\CF_n^-$.

Thus, it is enough to show the following lemma stating that the
forward likelihood ratios are uniformly tight in a strong sense. The estimates have to be uniform in all $g$-chains, since our substitutions $\mu=\nu\circ T^{-m}$ and $\tl\mu=\tl\nu\circ T^{-m}$ are valid only for a fixed $m$.

\begin{lem}\label{concl}
  Assume that $\var_n g$ is square summable and that $g$ is
  bounded away from zero. For all $\epsilon > 0$ there is a
  $K=K(\epsilon)$ such that
\begin{equation}\label{uniform}
\sup_{\tl\mu,\mu} \, \sup_n 
\tl\mu\left( \log Z_n(x; \tl\mu, \mu, \{\CF_n^-\}) > K \right) <
\epsilon, 
\end{equation}
where $\tl\mu$ and $\mu$ are chosen among all pairs of forward
$g$-chains. 
\end{lem}
\begin{proof}
We adapt the part of the proof in \cite{johob}, pp.\ 597--598. One could also adapt the more advanced theory of Shiryaev \cite{shiryaev} and his co-authors, as was done in \cite{jopp}.
  
  Given $\mu$ and $\tl\mu$ write (with notation from \cite{johob}) $M_n(x)=Z_n(x; \tl\mu, \mu,
  \{\CF_n^-\})$.
  We show that $\log M_n$ has a Doob decomposition (see e.g.\ \cite{williams}, pp.\ 120-121):  
  $$ \log M_n = A_n + \eta_n $$
  where $A_n$ is previsible with the uniform bound $A_n \leq C_1 \sum
  (\var_n g)^2$, and moreover, $\eta_n$ is a
  $\tl\mu$-martingale and uniformly bounded in $L_2(\tl\mu)$ with
  $\tl\mu(\eta_n^2) \leq C_2 \sum (\var_n g)^2$. 
It is the uniformity of the estimates which makes it possible to conclude the strong formulation of \eqref{uniform}. We define
$$P_n(x)=\frac{\mu[x_{-n},x_{-n+1},\ldots,x_{-1}]}{\mu[x_{-n+1},\ldots,x_{-1}]}$$
$$\tl P_n(x)=\frac{\tl\mu[x_{-n},x_{-n+1},\ldots,x_{-1}]}{\tl\mu[x_{-n+1},\ldots,x_{-1}]}$$
and note that $|\tl P_n -P_n|\leq \var_n g$.

The proof then follows exactly as in \cite{johob}, pp.\ 597--598.
\end{proof}

\section{Open questions}\noindent
In this section we present two questions which the authors find
interesting, as well as challenging, in the light of the present
investigation.
\subsection{Question 1} Does uniqueness of a $g$-measure imply
$$\mathcal\L^n f\to \int f \d \mu?$$ We assume that we have finitely
many symbols for the left shift map $T$ and a continuous and strictly
positive $g$-function.

In this paper we used that a unique measure in $K_g$ (Walters
notation), or equivalently, a unique $G$-measure (Brown and Dooley
\cite{brown}, Fan \cite{fan}, and others), implies that $\mathcal\L^n
f(x)\to \int f\d \mu$, uniformly in $x\in X^+$. It would be natural to
ask that if $K_g$ only contains the $g$-measure, then is this the
unique member also of $K_g$ (recall that a $g$-measure is always a
member of $K_g$)? Then it would follow that uniqueness of a
$g$-measure (under the assumptions given in our question) implies
convergence of the iterates of $\mathcal\L$, without imposing stronger
regularity on $g$ than continuity.

\subsection{Question 2} Under square summability of variations of the $g$-function, what is the rate of convergence of $\mathcal\L^n f(x)$ in the supremum-norm? Rates of convergence
in the case of {\em summable variations} are available, see 
Pollicott \cite{pollicott}.

{\bf Acknowledgements.} The authors are grateful for the hospitality shown by the University of Warwick and in particular for the conversations with Thomas Jordan, Mark Pollicott and Peter Walters.

\end{document}